\documentclass{gtart_h}
\def\ifplaintex{\expandafter\ifx\csname documentclass\endcsname\relax}

\def\gtp{{\mathsurround=0pt\it $\cal G\mskip-2mu$eometry \&\ 
$\cal T\!\!$opology $\cal P\!$ublications}}  % GT publications

\def\recd{{\small Received:\qua\receiveddate\ifx\reviseddate\relax
\else\qquad Revised:\qua\reviseddate\fi\par}} 

\def\lognumber#1{\def\thelognumber{#1}}
\def\volumenumber#1{\def\thevolumenumber{#1}}
\def\volumeyear#1{\def\thevolumeyear{#1}}
\def\papernumber#1{\def\thepapernumber{#1}}
\def\pagenumbers#1#2{\def\startpage{#1}\def\finishpage{#2}}
\def\published#1{\def\publishdate{#1}}

\def\received#1{\def\receiveddate{#1}}

\def\accepted#1{\def\accepteddate{#1}}
\def\asciititle#1{\def\theasciititle{#1}}

\long\def\asciiabstract#1{\long\def\theasciiabstract{#1}}
\def\asciikeywords#1{\def\theasciikeywords{#1}}

\let\\\par\let\thelognumber\relax\let\thevolumenumber\relax
\let\thepapernumber\relax\let\thevolumeyear\relax\let\startpage\relax
\let\finishpage\relax\let\publishdate\relax\let\receiveddate\relax
\let\reviseddate\relax\let\accepteddate\relax\let\theasciititle\relax
\let\theasciiauthors\relax
\let\theasciiabstract\relax\let\theasciikeywords\relax

\let\theasciiemail\relax

\ifplaintex
\font\logobig=cmssbx10 scaled 3836
\font\logomed=cmssbx10 scaled 2557
\else
\font\logobig=cmssbx10 scaled 4200
\font\logomed=cmssbx10 scaled 2800
\fi

\long\def\makeagttitle{   %%% start of definition of \makeagttitle
\count0=\startpage
\agt\hfill      %   Journal title (top left) 
\hbox to 45truept{\vbox to 0pt{\vglue -13truept{\logomed A\kern -.37em{\logobig 
T}\kern -.38em G}\vss}\hss}
\break
{\small Volume \thevolumenumber\ (\thevolumeyear)
\startpage--\finishpage\nl
Published: \publishdate}

\vglue .25truein

{\parskip=0pt\leftskip 0pt plus
1fil\def\\{\par\smallskip}{\Large\bf\thetitle}\par\medskip} \vglue
0.05truein

{\parskip=0pt\leftskip 0pt plus 1fil\def\\{\par}{\sc\theauthors}
\par\medskip}%
 
\vglue 0.03truein 

{\small\leftskip 25truept\rightskip 25truept{\bf Abstract}\stdspace\theabstract

{\bf AMS Classification}\stdspace\theprimaryclass
\ifx\thesecondaryclass\relax\else; \thesecondaryclass\fi\par
{\bf Keywords}\stdspace \thekeywords\par}\vglue 7truept

}   %%%% end of definition of \makeagttitle

\ifplaintex
\font\phead=cmsl9 scaled 950
\font\pnum=cmbx10 scaled 913
\font\pfoot=cmsl9 scaled 950
\headline{\vbox to 0pt{\vskip -4.5mm\line{\small\phead\ifnum
\count0=\startpage ISSN 1472-2739 (on-line) 1472-2747 (printed)
\hfill {\pnum\folio}\else\ifodd\count0\def\\{ }% 
\ifx\theshorttitle\relax\thetitle\else\theshorttitle\fi\hfill{\pnum\folio}
\else\def\\{ and }{\pnum\folio}\hfill\ifx\theshortauthors\relax\theauthors
\else\theshortauthors\fi\fi\fi}\vss}}
\footline{\vbox to 0pt{\vglue 0mm\line{\small\pfoot\ifnum\count0=\startpage
\copyright\ \gtp\hfill\else
\agt, Volume \thevolumenumber\ (\thevolumeyear)\hfill\fi}\vss}}
\else
\font=cmsl9 scaled 1050
\font\lnum=cmbx10 
\font=cmsl9 scaled 1050
\makeatletter
\def\@oddhead{{\small\lhead\ifnum\count0=\startpage ISSN 1472-2739 
(on-line) 1472-2747 (printed)\hfill {\lnum\number\count0}\else\ifodd\count0
\def\\{ }\ifx\theshorttitle\relax \thetitle \else\theshorttitle\fi\hfill
{\lnum\number\count0}\else\def\\{ and }{\lnum\number\count0}
\hfill\ifx\theshortauthors\relax 
\theauthors\else\theshortauthors\fi\fi\fi}}\def\@evenhead{\@oddhead}
\def\@oddfoot{\small\lfoot\ifnum\count0=\startpage\copyright\ \gtp\hfill\else
\agt, Volume \thevolumenumber\ (\thevolumeyear)\hfill\fi}
\def\@evenfoot{\@oddfoot}
\makeatother
\fi
\let\maketitlepage\makeagttitle
\let\makeshorttitle\maketitlepage
\let\maketitle\maketitlepage

\newwrite\gtoutfile
\long\gdef\makeheadfile{  %%% start of definition of \makeheadfile
{\def\\{, }\def\s{ }
\immediate\openout\gtoutfile head.xxx
\immediate\write\gtoutfile{Proxy-for: \ifx\theasciiauthors\relax
\theauthors\else\theasciiauthors\fi\s<\ifx\theasciiemail\relax\theemail\else\theasciiemail\fi>}
\immediate\write\gtoutfile{\noexpand\\}
\immediate\write\gtoutfile{Authors: \ifx\theasciiauthors\relax
\theauthors\else\theasciiauthors\fi}
{\def\\{ }\immediate\write\gtoutfile{Title: \ifx\theasciititle\relax
\thetitle\else\theasciititle\fi}}
\immediate\write\gtoutfile{Subj-class: GT or SG, GR etc}
\immediate\write\gtoutfile{MSC-class: \theprimaryclass\ifx\thesecondaryclass\relax\else, \thesecondaryclass\fi}
\immediate\write\gtoutfile{Journal-ref: Algebr. Geom. Topol. \thevolumenumber\s
(\thevolumeyear) \startpage-\finishpage}
\immediate\write\gtoutfile{Comments: Published by Algebraic and
Geometric Topology at}
\immediate\write\gtoutfile{\s\s\s  http://www.maths.warwick.ac.uk/agt/AGTVol\thevolumenumber/agt-\thevolumenumber-\thepapernumber.abs.html}
\immediate\write\gtoutfile{\noexpand\\}
\immediate\write\gtoutfile{}
\ifx\theasciiabstract\relax
\immediate\write\gtoutfile{\theabstract}\else
\immediate\write\gtoutfile{\theasciiabstract}\fi
\immediate\write\gtoutfile{}
\immediate\write\gtoutfile{\noexpand\\}
\immediate\write\gtoutfile{}
\immediate\closeout\gtoutfile}}  %%% end of definition of \makeheadfile

\def\maketitlepage{\makeagttitle\makeheadfile}
\let\makeshorttitle\maketitlepage
\let\maketitle\maketitlepage

\lognumber{47}
\volumenumber{4}
\volumeyear{2004}
\papernumber{47}
\pagenumbers{1103}{1109}
\received{4 August 2004} 
\accepted{16 November 2004}
\published{25 November 2004}

\usepackage{amssymb, amscd, amsmath, enumerate}

\newtheorem{theorem}{Theorem}
\newtheorem{lemma}[theorem]{Lemma}

\begin{document}

\title{An indecomposable $PD_3$-complex : II}
\asciititle{An indecomposable PD_3-complex : II}

\author{Jonathan A. Hillman }
\address{School of Mathematics and Statistics F07\\University of Sydney, 
NSW 2006, Australia}

\email{jonh@maths.usyd.edu.au}

\begin{abstract}
We show that there are two homotopy types of
$PD_3$-complexes with fundamental group $S_3*_{Z/2Z}S_3$,
and give explicit constructions for each,
which differ only in the attachment of the top cell.
\end{abstract}

\asciiabstract{We show that there are two homotopy types of
PD_3-complexes with fundamental group S_3*_{Z/2Z}S_3, and give
explicit constructions for each, which differ only in the attachment
of the top cell.}

\asciikeywords{Indecomposable, Poincare duality, PD_3-complex}

\maketitle

In \cite{[Hi05]} we showed that $\pi=S_3*_{Z/2Z}S_3$ satisfies the criterion
of \cite{[Tu90]} and thus is the fundamental group of a $PD_3$-complex.
As $\pi$ has infinitely many ends but is indecomposable,
this illustrates a divergence from the known properties of 3-manifolds,
and provides a counter-example to an old question of Wall \cite{[Wa67]}.
In particular, the Sphere Theorem does not extend to all $PD_3$-complexes.

Here we shall give an explicit description of a finite $PD_3$-complex
$Y$ realizing this group.
The construction is modelled on a similar construction for
a $PD_3$-complex $X$ with fundamental group $S_3$.
In each case the cellular chain complex of the universal
cover has the striking property that it is self-dual.
In \S2 we show a $PD_3$-complex with fundamental group $\pi$ must be orientable,
and we use Turaev's work to show there are two homotopy types of
such $PD_3$-complexes.
The 2-fold cover of $Y$ is homotopy equivalent to $L(3,1)\sharp{L(3,1)}$,
while a simple modification of our construction (suggested by the referee)
gives a $PD_3$-complex with 2-fold cover homotopy equivalent to
$L(3,1)\sharp{-L(3,1)}$.
(This group was first suggested as a test case in \cite{[Hi93]}.)

\section{A finite complex with group $S_3*_{Z/2Z}S_3$}

Let $G$ be a group and let $\Gamma=\mathbb{Z}[G]$,
$\varepsilon:C_1=\Gamma\to\mathbb{Z}$ and $I(G)=\mathrm{Ker}(\varepsilon)$
be the integral group ring,
the augmentation homomorphism and the augmentation ideal, respectively.
If $M$ is a left $\Gamma$-module $\overline{M}$ shall denote the conjugate
right module, with $G$-action given by $m.g=g^{-1}m$ for all $g\in G$ and
$m\in M$, and similarly $\overline{N}$ shall denote the conjugate
left module structure on a right $\Gamma$-module $N$.
If $C_*$ is a chain complex over $\Gamma$ with an
augmentation $\varepsilon:C_0\to\mathbb{Z}$ a {\it diagonal approximation}
is a chain homomorphism
$\Delta:C_*\to{C_*}\otimes_\mathbb{Z}C_*$ (with diagonal $G$-action)
such that $(\varepsilon\otimes1)\Delta=id_{C_*}=(1\otimes\varepsilon)\Delta$.

The cellular chain complex $C_*(\widetilde{K})$
for the universal covering space of a finite
2-complex $K$ determined by a presentation for a group
is isomorphic to the Fox-Lyndon complex of the presentation,
via an isomorphism carrying generators corresponding to based lifts of
cells of $K$ to the standard generators.

The symmetric group $S_3$ has a presentation
$\langle a,b\mid a^2,abab^{-2}\rangle$.
Let $\pi=S_3*_{Z/2Z}S_3$, with presentation $\langle a,b,c\mid r,s,t\rangle$,
where $r=a^2$, $s=abab^{-2}$ and $t=acac^{-2}$.
The two obvious embeddings of $S_3$ into $\pi$ admit retractions,
as $\pi/\langle\langle{b}\rangle\rangle\cong
\pi/\langle\langle{c}\rangle\rangle\cong{S_3}$.
Let $A$, $B$ and $C$ be the cyclic subgroups generated by the images of
$a$, $b$ and $c$, respectively.
The inclusions of $A$ into $S_3$ and $\pi$ induce isomorphisms on
abelianization, while the commutator subgroups are $S_3'=B$ and $\pi'=B*C$.
Thus these groups are semidirect products: $S_3\cong B\rtimes(Z/2Z)$ and
$\pi\cong(B*C)\rtimes{Z/2Z}$.
In particular, $\pi$ is virtually free, and so has infinitely many ends.
However it follows easily from the Grushko-Neumann Theorem that $\pi$ is
indecomposable. (See \cite{[Hi05]}).

The above presentations determine finite 2-complexes $K$ and $L$,
with fundamental groups $S_3$ and $\pi$, respectively.
There are two obvious embeddings of $K$ as a retract in $L$,
with retractions $r_b,r_c:L\to{K}$ given by collapsing the pair of cells
$\{c,t\}$ and $\{b,s\}$, respectively.

The chain complex $C_*(\widetilde{K})$ has the form
\begin{equation*}
\begin{CD}
\mathbb{Z}[S_3]^2@>\partial_2>>\mathbb{Z}[S_3]^2@>\partial_1>>\mathbb{Z}[S_3],
\end{CD}
\end{equation*}
where $\partial_1(1,0)=a-1$, $\partial_1(0,1)=b-1$,
$\partial_2(1,0)=(a+1,0)$ and  $\partial_2(0,1)=(b^2a+1,a-b-1)$.
The 2-chain $\psi=(a-1,-ba+a+b^2-b)$ is a 2-cycle,
and so determines an element of $\pi_2(K)=H_2(\widetilde{K};\mathbb{Z})$,
by the Hurewicz Theorem.
Let $X=K\cup_\psi{e^3}$, and let $C_*$ be the cellular chain complex
for the universal cover $\widetilde{X}$.
(Thus $C_i=C_i(\widetilde{K})$ for $i\leq2$ and $C_3\cong\mathbb{Z}[S_3]$).
The dual cochain complex ${C}^*={Hom_\Gamma(C_*,\mathbb{Z}[S_3])}$ is a complex
of right $\mathbb{Z}[S_3]$-modules.

We shall define new bases which display the structure of $C_*$ to
better advantage, as follows. Let $e_1=(1,0)$ and
$e_2=(-ba-b^2,1)$ in $C_1$ and $f_1=(1,0)$ and $f_2=(0,-a)$ in
$C_2$, and let $g$ be the generator of $C_3$ corresponding to the
top cell. Then $\partial_1{e_1}=a-1$,
$\partial_1{e_2}=-b^2a+ba+b^2-1$, $\partial_2{f_1}=(a+1)e_1$,
$\partial_2{f_2}=(b^2a+a-1)e_2$, and
$\partial_3g=\psi=(a-1)f_1+(-b^2a+ba+b-1)f_2$. The matrix for
$\partial_2$ with respect to the bases $\{\tilde e_i\}$ and
$\{\tilde f_j\}$ is diagonal, and is hermitian with respect to the
canonical involution of $\mathbb{Z}[S_3]$, while the matrix for
$\partial_3$ is the conjugate transpose for that of $\partial_1$.
Hence the chain complex $\overline{C^{3-*}}$ obtained by
conjugating and reindexing the cochain complex $C^*$ is isomorphic
to $C_*$.

Let $\beta=b^2+b+1$ and $\nu=\Sigma_{s\in{S_3}}s=\beta(a+1)$.

\begin{lemma}
The complex $X$ is a $PD_3$-complex with $\widetilde{X}\simeq{S^3}$.
\end{lemma}

\begin{proof}
Since $C_*$ is the cellular chain complex of a 1-connected
cell complex $H_0(C_*)\cong\mathbb{Z}$ and $H_1(C_*)=0$.
If $\partial_2(rf_1+sf_2)=0$ then $r(a+1)=0$ and $s(b^2a+a-1)=0$.
Now the left annihilator ideals of $a+1$ and $b^2a+a-1$ in $\mathbb{Z}[S_3]$
are principal left ideals, generated by $a-1$ and
$(b-1)(ba-1)$, respectively.
Hence $r=p(a-1)$ and $s=q(b-1)(ba-1)$ for some $p,q\in\mathbb{Z}[B]$.
A simple calculation gives $\partial_3((p(ba+b+1)+q(ba+b))g)=rf_1+sf_2$
and so $H_2(C_*)=0$.

If $\partial_3hg=0$ then $h(a-1)=0$, so $h=h_1(a+1)$ for some
$h_1\in\mathbb{Z}[B]$, and $h(b^2a-ba-b+1)=0$.
Now $h(b^2a-ba-b+1)=h_1(1-b)(a+b+1)$, so $h_1(1-b)=0$.
Therefore $h_1=m\beta$ for some $m\in\mathbb{Z}$,
so $h=m\nu$ and $H_3(C_*)=\mathbb{Z}[S_3]\nu{g}\cong\mathbb{Z}$.
Hence $\widetilde{X}\simeq{S^3}$.
Now $H_3(X;\mathbb{Z})=H_3(\mathbb{Z}\otimes_{\mathbb{Z}[S_3]}C_*)=
\mathbb{Z}[1\otimes{g}]$ and $tr([1\otimes{g}])=\nu{g}$,
where $tr:H_3(X;\mathbb{Z})\to H_3(\widetilde{X};\mathbb{Z})$ is
the transfer homomorphism.
The homomorphisms from $H^q(\overline{C^*})$ to $H_{3-q}(C_*)$
determined by cap product with $[X]=[1\otimes{g}]$
may be identified with the Poincar\'e duality isomorphisms for
$\widetilde{X}$, and so $X$ is a $PD_3$-complex.
\end{proof}

The verification that $\widetilde{X}\simeq S^3$ is essentially due to
\cite{[Sw60]} and the fact that $X$ is a $PD_3$-complex is due to \cite{[Wa67]}.
The only novelty here is the diagonalization of $\partial_2$,
which was a guiding feature in the study of $\pi=S_3*_{Z/2Z}S_3$.

Let $\Pi=\mathbb{Z}[\pi]$.
The cellular chain complex for the universal covering space
$\widetilde{L}$ has the form
\begin{equation*}
\begin{CD}
\Pi^3@>\partial_2>>\Pi^3@>\partial_1>>\Pi.
\end{CD}
\end{equation*}
The differentials are given by
$\partial_1(1,0,0)=a-1$, $\partial_1(0,1,0)=b-1$ and
$\partial_1(0,0,1)=c-1$,
$\partial_2(1,0,0)=(a+1,0,0)$, $\partial_2(0,1,0)=(b^2a+1,a-b-1,0)$
and $\partial_2(0,0,1)=(c^2a+1,0,a-c-1)$.
In particular, $H_2(\widetilde{L};\mathbb{Z})=\mathrm{Ker}(\partial_2)$.

Let $\theta=(a-1,-ba+a+b^2-b,-ca+a+c^2-c)$.
Then $\partial_2(\theta)=0$,
and so $\theta$ determines an element of $\pi_2(L)=H_2(\widetilde{L};\mathbb{Z})$,
by the Hurewicz Theorem.
Let $Y=L\cup_\theta{e^3}$ and let $D_*$ be the cellular chain complex
for the universal covering space $\widetilde{Y}$.

\begin{tabular}{lll}
$\phantom{99999}$ & Let & Then\\
&$\tilde e_1=(1,0,0)$ & $\partial_1{\tilde e_1}=a-1$\\
&$\tilde e_2=(-ba-b^2,1,0)\phantom{99}$ & $\partial_1{\tilde e_2}={ba-b^2a+b^2-1}$\\
&$\tilde e_3=(-ca-c^2,0,1)$ & $\partial_1{\tilde e_3}={ca-c^2a+c^2-1}$\\
&$\tilde f_1=(1,0,0)$ & $\partial_2{\tilde f_1}=(a+1)\tilde e_1$\\
&$\tilde f_2=(0,-a,0)$ & $\partial_2{\tilde f_2}=(b^2a+a-1)\tilde e_2$\\
&$\tilde f_3=(0,0,-a).$ & $\partial_2{\tilde f_3}=(c^2a+a-1)\tilde e_3.$\\
\end{tabular}

\noindent  Moreover $\theta=
(a-1)\tilde f_1+(-b^2a+ba+b-1)\tilde f_2+(-c^2a+ca+c-1)\tilde f_3.$
Let ${D}^*={Hom_\Gamma(D_*,\Pi)}$ be the cochain complex dual to $D_*$.
Then it is easily seen that $\overline{D^{*}}\cong D_{3-*}$.

\begin{theorem}
 The complex $Y$ is a $PD_3$-complex.
\end{theorem}

\begin{proof}
Clearly $H_0(D_*)\cong\mathbb{Z}$ and $H_1(D_*)=0$.
The argument of the first part of Lemma 1 extends immediately
to show that the kernel of $\partial_2$ is generated by $(a-1)\tilde f_1$,
$(b-1)(ba-1)\tilde f_2$ and $(c-1)(ca-1)\tilde f_3$.
Hence these elements represent generators for $H_2(D_*)$.
Let $\tilde g$ be the generator for $D_3$ corresponding to the top cell,
so that $\partial_3\tilde g=\theta$.
Note that the image of $g$ in $\mathbb{Z}\otimes_\varepsilon D_3$ is a cycle, and
represents a generator for $H_3(Y;\mathbb{Z})=H_3(\mathbb{Z}\otimes_\varepsilon D_*)$.
If $h\theta=0$ then (as in Lemma 1) $h=h_1(a+1)$ for some $h_1\in\mathbb{Z}[B*C]$
such that $h_1(b-1)=h_1(c-1)=0$.
It follows that $h_1=0$.
Hence $\partial_3$ is injective and so $H_3(D_*)=0$.

Let $\hat1$, $\hat e^*$, $\hat f^*$ and $\hat g$ denote the bases of
${D}^*$ dual to the above bases for $D_*$.
Let $\Delta$ be a diagonal approximation for $D_*$ and suppose that
$\Delta(\tilde{g})=\Sigma_{0\leq{q}\leq3}\Sigma_{i\in I(q)}x_i\otimes{y}_i$,
where $x_i\in D_q$ and $y_i\in D_{3-q}$, for all $i\in I(q)$ and $0\leq q\leq3$.
Then $\Sigma_{i\in I(3)}x_i=\tilde{g}$.
Let $r_i=\hat{g}(x_i)$ for $i\in I(3)$ and let $\tilde\xi$ denote the image
of $\tilde g$ in $H_3(Y;\mathbb{Z})=\mathbb{Z}\otimes_\varepsilon D_3$.
Then $\varepsilon(\hat{g}\cap\tilde\xi)=
\varepsilon(\Sigma_{i\in I(3)}\overline{r_i}y_i)=
\varepsilon(\Sigma_{i\in I(3)}\overline{r_i})=
\varepsilon(\overline{\hat{g}(\tilde{g})})=1$,
and so $\hat{g}\cap\tilde\xi$ generates $H_0(D_*)$.
Since $H_1(D_*)=H_3(D_*)=H^0(\overline{D^*})=H^2(\overline{D^*})=0$,
$-\cap\tilde\xi$ induces isomorphisms
$H^q(\overline{D^*})\cong H_{3-q}(D_*)$ for all $q\not=1$.
The remaining case follows as in \cite{[Tu90]} from the facts that
$\overline{D^{*}}\cong D_{3-*}$ and $\Delta$ is chain homotopic to $\tau\Delta$,
where $\tau:D_*\otimes{D_*}\to D_*\otimes{D_*}$ is the transposition defined by
$\tau(\alpha\otimes\omega)=(-1)^{pq}\omega\otimes\alpha$ for all
$\alpha\in D_p$ and $\omega\in D_q$.
Thus $Y$ is a $PD_3$-complex.
\end{proof}

Can the last step of this argument be made more explicit?
The work of Handel \cite{[Ha93]} on diagonal approximations for dihedral
groups may be adapted to give the following formulae for a diagonal
approximation for the truncation to degrees $\leq2$ of $D_*$
which is compatible with the above two embeddings of $K$ as a retract in $L$:

\leftskip 25pt

\noindent $\Delta(1)=1\otimes 1$

\noindent $\Delta(\tilde e_1)=\tilde e_1\otimes{a}+1\otimes{\tilde e_1},$

\noindent
$\Delta(\tilde{e_2})=\tilde{e_2}\otimes1-ba\tilde{e_1}\otimes(b-1)
-b^2\tilde{e_1}\otimes(b^2a-1)-(ba-b)\otimes{ba\tilde{e_1}}$

$\phantom{999}-(b^2-b)\otimes{b^2\tilde{e_1}}+
b\otimes\tilde{e_2},$

\noindent
$\Delta(\tilde{e_3})=\tilde{e_3}\otimes1-ca\tilde{e_1}\otimes(c-1)
-c^2\tilde{e_1}\otimes(c^2a-1)-(ca-c)\otimes{ca\tilde{e_1}}$

$\phantom{999}-(c^2-c)\otimes{c^2\tilde{e_1}}+c\otimes\tilde{e_3},$

\noindent $\Delta(\tilde f_1)=\tilde f_1\otimes1+
\tilde e_1\otimes{a\tilde e_1}+1\otimes{\tilde f_1}$,

\noindent $\Delta(\tilde{f}_2)=\tilde{f}_2\otimes{a}
+(b^2+b)\tilde{f}_1\otimes{(a-ba)}+(b^2a+b^2)\tilde{f}_2\otimes{(a-ba)}$

$\phantom{999}+((ba+b^2-1)\tilde{e}_1+\tilde{e}_2)\otimes{((b^2a)\tilde{e}_1
+ba\tilde{e}_2)}$

$\phantom{999}-((b^2a+1)\tilde{e}_1+ba\tilde{e}_2)\otimes{((ba+a+b^2+b)\tilde{e}_1+(b^2a+a)\tilde{e}_2)}$

$\phantom{999}-((a+b)\tilde{e}_1+b^2a\tilde{e}_2)\otimes{((ba+b^2)\tilde{e}_1+a\tilde{e}_2)}-(a+1)\tilde{e}_1\otimes{\tilde{e}_1}$

$\phantom{999}+(a-b)\otimes(b^2+b)\tilde{f}_1+(a-b)\otimes(b^2a+b^2)\tilde{f}_2
+a\otimes{\tilde{f}_2}$ \quad and

\noindent $\Delta(\tilde{f}_3)=\tilde{f}_3\otimes{a}
+(c^2+c)\tilde{f}_1\otimes{(a-ca)}+(c^2a+c^2)\tilde{f}_3\otimes{(a-ca)}$

$\phantom{999}+((ca+c^2-1)\tilde{e}_1+\tilde{e}_3)\otimes{((c^2a)\tilde{e}_1+ca\tilde{e}_3)}$

$\phantom{999}-((c^2a+1)\tilde{e}_1+ca\tilde{e}_3)\otimes{((ca+a+c^2+c)\tilde{e}_1+(c^2a+a)\tilde{e}_3)}$

$\phantom{999}-((a+c)\tilde{e}_1+c^2a\tilde{e}_3)\otimes{((ca+c^2)\tilde{e}_1+a\tilde{e}_3)}
-(a+1)\tilde{e}_1\otimes{\tilde{e}_1}$

$\phantom{999}+(a-c)\otimes(c^2+c)\tilde{f}_1+(a-c)\otimes(c^2a+c^2)\tilde{f}_3
+a\otimes{\tilde{f}_3}$

\smallskip\leftskip 0pt

These formulae were derived from the work of Handel by using the canonical
involution of $\mathbb{Z}[S_3]$ to switch right and left module structures
and showing that $C_*$ is a direct summand of a truncation of the Wall-Hamada
resolution for $S_3$.
(In Handel's notation $a=y$, $b=x$, $e_1=c_1^2$, $e_2=-c_1^1-c_1^2(x+xy)$,
$f_1=c_2^3$,
$f_2=-c_2^1y+c_2^2x^2-c_2^3y$ and $g=-(c_3^1+c_3^3)(x+y)-c_3^4y$).
Handel's work also leads to a formula for $\Delta(g)$,
but it is not clear what $\Delta(\tilde{g})$ should be.

\section{Other $PD_3$-complexes with this group}

Having constructed one $PD_3$-complex with group $\pi$
one may ask how many there are.
Any such $PD_3$-complex must be orientable.
For let $w_1:\pi\to\{\pm1\}$ be a homomorphism and define an involution on
$\Gamma$ by $\bar g=w_1(g)g^{-1}$, for all $g\in\pi$.
Let $w=w_1(a)$ and $R=\mathbb{Z}[\pi/\pi']=\mathbb{Z}[a]/(a^2-1)$.
Let $J=\mathrm{Coker}(\overline{\partial_2}^{tr})$,
where $\partial_2:\Pi^3\to\Pi^3$ is the presentation matrix
for $I(\pi)$ given in \S1.
Then $R\otimes_\Gamma I(\pi)\cong R/(a+1)\oplus(R/(a+1,3))^2$,
while $R\otimes_\Gamma J\cong R/(a+w)\oplus(R/(a+w,3))^2$.
If the pair $(\pi,w_1)$ is realized by a $PD_3$-complex then $I(\pi)$ and $J$
are projective homotopy equivalent \cite{[Tu90]}.
But then $R\otimes_\Gamma I(\pi)$ and $R\otimes_\Gamma{J}$
are projective homotopy equivalent $R$-modules, and so we must have $w=1$.

If $W$ is an oriented $PD_3$-complex with fundamental group $G$ and
$c_W:W\to K(G,1)$ is a classifying map let $\mu(W)=c_{W*}[W]\in H_3(W;\mathbb{Z})$.
Two such $PD_3$-complexes $W_1$ and $W_2$ are homotopy equivalent if and only
$\mu(W_1)$ and $\mu(W_2)$ agree up to sign and the action of $Out(G)$.
Turaev constructed an isomorphism $\nu_C$ from $H_3(G;\mathbb{Z})$ to a group
$[F^2(C),I(G)]$ of projective homotopy classes of module homomorphisms
and showed that $\mu\in H_3(G;\mathbb{Z})$ is the image of the orientation class
of a $PD_3$-complex if and only if $\nu_C(\mu)$ is the class of a
self-homotopy equivalence \cite{[Tu90]}.

When $G=\pi=S_3*_{Z/2Z}S_3$ we have $F^2(C)\cong I(\pi)$,
and $H_3(\pi;\mathbb{Z})\cong{H_3(\pi';\mathbb{Z})}\oplus{H_3(Z/2Z;\mathbb{Z})}
\cong(Z/3Z)^2\oplus(Z/2Z)$.
Let $W'$ be the double cover of $W$, with fundamental group
$\pi'\cong(Z/3Z)*(Z/3Z)$.
Then $W'$ is a connected sum, by Theorem 1 of \cite{[Tu90]},
and so it is homotopy equivalent to one of the 3-manifolds
$L(3,1)\sharp{L(3,1)}$ and $L(3,1)\sharp-L(3,1)$.
(These may be distinguished by the torsion linking forms on
their first homology groups).
In particular, $\mu(W')$ has nonzero entries in each summand.
Since $\mu(W')$ is the image of $\mu(W)$ under the transfer
to $H_3(\pi';\mathbb{Z})\cong(Z/3Z)^2$ the image of $\mu(W)$ in
each $Z/3Z$-summand must be nonzero.
Let $u\in H^1(W;\mathbb{F}_2)$ correspond to the abelianization homomorphism.
Since $\beta_2(W;\mathbb{F}_2)=\beta_1(W;\mathbb{F}_2)=1=\beta_2(\pi;\mathbb{F}_2)$
we have $u^2\not=0$, and so $u^3\not=0$, by Poincar\'e duality.
It follows easily that the image of $\mu(W)$ in the
$Z/2Z$-summand must be nonzero also.
(Note that $W'$ is $\mathbb{Z}_{(2)}$-homology equivalent to $S^3$
and so $W$ is $\mathbb{Z}_{(2)}[Z/2Z]$-homology equivalent to $RP^3$).
Since reversing the orientation of $W$ reverses that of $W'$,
we may conclude that there are at most two distinct homotopy types
of $PD_3$-complexes with fundamental group $\pi$,
and that they may be detected by their double covers.

The retractions $r_b$ and $r_c$ of $L$ onto $K$ extend to maps
$r_b, r_c:Y\to X$.
These maps induce the same isomorphism $H_3(Y;\mathbb{Z})\cong H_3(X;\mathbb{Z})$,
and so their lifts to the double covers induce the same
isomorphism $H_3(Y';\mathbb{Z})\to H_3(X';\mathbb{Z})$.
Hence $ Y'\simeq {L(3,1)}\sharp{L(3,1)}$, rather than $L(3,1)\sharp{-L(3,1)}$.
The referee has pointed out that if we use
$\xi=(a-1)\tilde f_1+(-b^2a+ba+b-1)\tilde f_2-(-c^2a+ca+c-1)\tilde f_3$
instead of $\theta$ (changing only the sign of the final term) then
$Z=L\cup_\xi{e^3}$ is another $PD_3$-complex with $\pi_1(Z)\cong\pi$,
and a similar argument shows that the double cover is now
$Z'\simeq{L(3,1)}\sharp{-L(3,1)}$.

The question of whether every aspherical $PD_3$-complex is homotopy equivalent
to a 3-manifold remains open.
The recent article \cite{[Wa04]} gives a comprehensive survey of Poincar\'e
duality in dimension 3,
emphasizing the role of the JSJ decomposition in relation to this question.

\medskip
{\bf Acknowledgement}\qua
I would like to thank the referee for suggesting the modification
giving the example covered by ${L(3,1)}\sharp{-L(3,1)}$,
and for other improvements to the exposition.

\Addresses\recd

\enddocument